\documentclass [11pt,titlepage]{article}
\usepackage{amsmath}
\usepackage{amssymb}
\begin{document}

\begin{center}
\begin{Large} 
{\bf FIRST CASE OF FERMAT'S LAST THEOREM}
\end{Large}

\addvspace{0.5cm}
{\bf Joseph Amal Nathan}

Reactor Physics Design Division, Bhabha Atomic Research Centre, Mumbai-400085, India\\
{\bf email:}josephan@magnum.barc.ernet.in
\end{center}

\noindent
{\bf Abstract:} In this paper two conjectures are proposed based on which we can prove the first case of Fermat's Last Theorem(FLT) for all primes $p \equiv -1 (\bmod~6)$. With Pollaczek's result {\bf [1]} and the conjectures the first case of FLT can be proved for all primes greater than $3$. With a computer Conjecture1 was verified to be true for primes $\leq 2437$ and Conjecture2 for primes $\leq 100003$. 

\noindent
{\bf Fermat's Last Theorem(FLT):} {\it The equation $x^l+y^l=z^l$ with integral $l > 2$, has no solution in positive integers} $x,y,z$.

\noindent
There is no loss in generality if $x,y,z$ are pairwise prime and, if $l=4$ and all odd primes only. Let $p$ denote a prime greater than $3$. The results presented here are when $p \nmid xyz$, known as `~{\it first case of FLT}~'. We will denote greatest common divisor of integers $\alpha,\beta,\gamma,...$ by ($\alpha,\beta,\gamma,$...).

\noindent
Let us define a {\it Polynomial} $F_n(x,y) = (x+y)^n-x^n-y^n$ for any positive odd integer $n$. Expanding $(x+y)^n$ and collecting terms having same binomial coefficient we get,
\begin{equation}
F_n(x,y)=x y (x+y) \left[ \sum_{i=1}^{(n-1)/2} {^nC_i}~ x^{i-1} y^{i-1} \left( \sum_{j=0}^{n-1-2i}(-1)^j x^{n-2i-j-1} y^j \right) \right].
\end{equation}
We will show $F_n(x,y)$ can be written in terms of $x,y$ and $(x+y)$,
\begin{equation}
F_n(x,y)=x y (x+y)\sum_{i=1}^{(n-1)/2} W_{i-1}~ x^{i-1} y^{i-1} (x+y)^{n-2i-1},
\end{equation}
\noindent
where, $W_j$, for $j=0,1,...,\frac {n-3}{2}$, are integers. Comparing (1) and (2) we get, 
\begin{equation}
W_j={^nC_{j+1}} - \sum_{k=1}^{j} \left[ {(-1)^{j+k}}~ {^nC_k} + {^{n-2j+2k-3}C_k}~ W_{j-k} \right],
\end{equation}
\noindent
which shows all $W_j$'s are integers. We see $W_0=n$ and evaluating the last coefficient,
$$\left. \frac {F_n(x,y)}{x+y} \right|_{x=-y} = (-1)^{\frac{n-1}{2}}~
W_{\frac{n-3}{2}}~ x^{n-1} = \left[(x+y)^{n-1} - \left. \left( \frac {x^n+y^n} {x+y} \right) \right] \right|_{x=-y} = -n x^{n-1},$$
\noindent
we get
\begin{equation}
W_{\frac{n-3}{2}} = (-1)^{\frac{n-3}{2}}~ n.
\end{equation}

\noindent
{\bf Lemma1:} {\it Given nonzero integers $x,y$ such that $(x,y)=1$ for any positive odd integer $n$},
$$\left( x+y, \frac {x^n+y^n}{x+y} \right) = (x+y,n).$$ 
\noindent
Let $(x+y,n)=g$.

\noindent
{\it Proof 1}.~From (2) and (4), $\frac{F_n(x,y)}{g(x+y)} =\frac{(x+y)}{g}\left[ \sum_{i=1}^{(n-3)/2} W_{i-1}~ x^i y^i (x+y)^{n-2i-2}\right]-\frac{n}{g}(-xy)^{\frac{n-1}{2}}$, any divisor of $\frac{(x+y)}{g}$ will divide all terms except last, giving ${\left(\frac{x+y}{g},\frac{F_n(x,y)}{g(x+y)}\right) = 1}$. So $\left( \frac{x+y}{g}, \frac {x^n+y^n}{g(x+y)} \right) = 1.$

\noindent
{\it Proof 2}.~Since, $\frac {x^n+y^n}{g(x+y)} = \frac{(x+y)}{g} \left[\sum_{i=1}^{n-1}(-1)^{i-1}~i~x^{n-i-1} y^{i-1} \right]+\frac{n}{g}y^{n-1},$ any divisor of $\frac{(x+y)}{g}$ will divide all terms except last. Hence $\left( \frac{x+y}{g}, \frac {x^n+y^n}{g(x+y)} \right) = 1.~~~~~~~~~~~~~~~~~~~~~~~~~~~~~~~~~~~~~~~~~~~~~~~~~~~~\Box$

\noindent
Now let $n=p$, then from (3), $w_j=W_j/p$, for $j=0,1,...,\frac {n-3}{2}$, are integers. Define {\it Polynomials} $G_p(x,y)= \frac {F_p(x,y)}{p x y(x+y)}$, and $f_{x,y}=x^2+xy+y^2$. We will now show, {\it $f^\epsilon_{x,y}$ is a factor of $G_p(x,y)$ for any prime $p>3$, when exponent $\epsilon=1$}. We will also show $\epsilon=2$ if and only if $p \equiv 1 (\bmod~6)$ and $\epsilon \leq 2$. Solving for $x$ in $f_{x,y}=0$, gives $x=-y e^{\pm \frac {i \pi}{3}}$. Now to show $f^\epsilon_{x,y} \mid G_p(x,y)$ it is sufficient to show $\left. \frac {F_p(x,y)}{f^{\epsilon-1}_{x,y}} \right|_{x =-y e^{\pm \frac{i \pi}{3}}}=0$. When $\epsilon=1$ subsituting for $x$ in $F_p(x,y) = (x+y)^n-x^n-y^n$, we get 
$$\left. F_p(x,y) \right|_{x =-y e^{\pm \frac{i \pi}{3}}}=y^p \left[ 2 \cos \frac {p \pi}{3}-1 \right],$$
\noindent
which is zero for any prime $p \equiv \pm 1(\bmod~6)$. Since $F_p(x,y)/f ^{\epsilon-1}_{x,y}$ at $x =-y e^{\pm \frac{i \pi}{3}}$ is indeterminant for $\epsilon > 1$, applying L'Hospital's rule when $\epsilon=2$ we get,
$$\left. \frac {F_p(x,y)}{f_{x,y}} \right|_{x=-y e^{\pm \frac{i \pi}{3}}}= \left. \frac {[F_p(x,y)]^{\prime}}{[f_{x,y}]^{\prime}} \right|_{x=-y e^{\pm \frac{i \pi}{3}}}= \mp 2 i p y^{p-2} \left[\frac {\sin \frac {(p-1) \pi}{3}}{1-2~e^{\pm \frac {i \pi}{3}}} \right],$$
\noindent
which is zero if and only if $p \equiv 1 (\bmod~6)$. When $\epsilon=3$ we get,
$$\left. \frac {F_p(x,y)}{f^2_{x,y}} \right|_{x=-y e^{\pm \frac{i \pi}{3}}}= \left. \frac {[F_p(x,y)]^{\prime \prime}}{[f^2_{x,y}]^{\prime \prime}} \right|_{x=-y e^{\pm \frac{i \pi}{3}}}= p (p-1) y^{p-4} \left[\frac {\cos \frac {(p-2) \pi}{3}}{ \left( 1-2~e^{\pm \frac {i \pi}{3}} \right)^2 } \right],$$
\noindent
which is not zero for $p \equiv \pm 1(\bmod~6)$. So $\epsilon \leq 2$ and $\epsilon = 1$ or $2$ for primes $p \equiv -1$ or $+1 (\bmod~6)$ respectively.    

\noindent
{\bf Lemma2:} {\it If the first case of FLT is false, then $p^2 \mid G_p(x,y)$ for primes $p \geq 3$}. 

\noindent
Since first case of FLT is false we have nonzero pairwise prime integers $x,y,z$ satisfying equation $x^p+y^p=z^p$. Define $m=x+y-z$, from Fermat's little theorem $p \mid m$. Let $p^ \lambda \mid m$ first we will show $\lambda \geq 2$. Since $p \nmid xyz$ we have $(x+y,p)=1$ and  from Lemma1
\begin{equation}
(x+y)=a^p;~~~~(z-x)=b^p;~~~~(z-y)=c^p, 
\end{equation}
\begin{equation}
\frac{x^p+y^p}{x+y}=r^p;~~~~\frac{z^p-x^p}{z-x}=s^p;~~~~\frac{z^p-y^p}{z-y}=t^p,  
\end{equation}
\begin{equation}
z=ar;~~~~~~y=bs;~~~~~~x=ct,
\end{equation}
\noindent
where $a, b, c, r, s, t$ are positive integers $>1$ and $p \nmid a, b, c, r, s, t$. Substituting (5) and (7) in $m=x+y-z$ we get, $m=a(a^{p-1}-r)=b(s-b^{p-1})=c(t-c^{p-1})$. Since $p \mid m$, we have $(a^{p-1})^p \equiv r^p(\bmod~p^2) \Rightarrow r^p \equiv 1(\bmod~p^2)$ and substituting for $r^p$ from (6) we get,
\begin{equation}
z^p \equiv (x+y)(\bmod~p^2).
\end{equation}
\noindent
Similarly 
\begin{equation}
y^p \equiv (z-x)(\bmod~p^2),
\end{equation}
\begin{equation}
x^p \equiv (z-y)(\bmod~p^2).
\end{equation}
\noindent
Taking (8)-(9)-(10), we get
\begin{equation}
m \equiv 0(\bmod~p^2),
\end{equation}
\noindent
so $\lambda \geq 2$. From (8)+(9)+(10), (8)+(9)-(10) and (8)-(9)+(10) we get $2z^p \equiv 2z(\bmod~p^2),~2y^p \equiv 2y(\bmod~p^2)$ and $2x^p \equiv 2x(\bmod~p^2)$ respectively, which gives
\begin{equation}
z^{p-1} \equiv 1(\bmod~p^2);~~~~y^{p-1} \equiv 1 (\bmod~p^2);~~~~x^{p-1} \equiv 1(\bmod~p^2).
\end{equation}
Now from (11) we have $(x+y)^p-z^p \equiv 0(\bmod~p^3) \Rightarrow p^2 \mid xy(x+y)G_p(x,y)$. Since $p \nmid xyz \Rightarrow p^2 \mid G_p(x,y)~~~~~~~~~~~~~~~~~~~~~~~~~~~~~~~~~~~~~~~~~~~~~~~~~~~~~~~~~~~~~~~~~~~~~~~~~~~~~~~~~~~~~~~~~~~~~~~~~~~~~~~~~~~\Box$

\noindent
{\bf Corollary1:} {\it For prime $3$ the first case of FLT is true}

\noindent
{\it Proof}.~Assume the first case of FLT is false for $p=3$. We have $G_3(x,y)=\frac {(x+y)^3-x^3-y^3}{3xy(x+y)}=1$. But from Lemma2, $3^2 \mid G_3(x,y)$, which is not possible. $~~~~~~~~~~~~~~~~~~~~~~~~~~~~~~~~~~~~~~~~~~~~~~~~~~~~~~~~~~~~~~~~~~~~~~~~~~~~~~~~~~~~~~~\Box$

\noindent
{\bf Lemma3:} {\it For any nonzero positive coprime integers $x,y$ there is no prime $\xi \equiv -1(\bmod~6)$ such that $\xi \mid (x^2+xy+y^2)$}.
 
\noindent
{\it Proof}.~We have $\xi +1=6 \tau$, where $\tau$ is a positive integer $\geq 1$. Let $x+y=\eta$. Assume $f_{x,y} = \eta^2-xy \equiv 0(\bmod~\xi)$. Since $(x,y)=1$ we get $(x,\xi)=(y,\xi)=(\eta,\xi)=(x,\eta)=(y,\eta)=1$ and the following congruences,
\begin{equation}
x^2 \equiv -y\eta(\bmod~\xi);~~~y^2 \equiv -\eta x(\bmod~\xi);~~~\eta^2 \equiv xy(\bmod~\xi),
\end{equation}
\begin{equation}
x^2+y^2+\eta^2 \equiv 0(\bmod~\xi). 
\end{equation}
\noindent
From (13) we get 
\begin{equation}
x^3 \equiv y^3 \equiv -\eta^3 (\bmod~\xi).
\end{equation}
\noindent
Using Fermat's little theorem and (15) in (14) we have,
$$x^2+y^2+\eta^2 \equiv x^{\xi+1}+y^{\xi+1}+\eta^{\xi+1} \equiv (x^{3\tau})^2 +(y^{3\tau})^2+(\eta^{3\tau})^2 \equiv 3x^{6\tau} \equiv 0(\bmod~\xi),$$
\noindent
which is contradictary to $\xi \nmid xy \eta.~~~~~~~~~~~~~~~~~~~~~~~~~~~~~~~~~~~~~~~~~~~~~~~~~~~~~~~~~~~~~~~~~~~~~~~~~~~~~~~\Box$

\noindent
{\bf Corollary2:} {\it For prime $5$ the first case of FLT is true}

\noindent
{\it Proof}.~Assume the first case of FLT is false for $p=5$. We have $G_5(x,y)=x^2+xy+y^2$ and $5 \equiv -1 (\bmod~6)$. But from Lemma2 we have $5^2 \mid G_5(x,y)$, contradicting Lemma3. $~~~~~~~~~~~~~~~~~~~~\Box$

\noindent
Define $H_p(x,y)=G_p(x,y)/f^\epsilon_{x,y}=\frac {F_p(x,y)}{pxy(x+y)f^\epsilon_{x,y}}$. It is obvious that $H_p(x,y)$ can be defined only for primes $p>3$. Now we propose the following {\it conjecture}.
 
\noindent
{\bf Conjecture1:} {\it For nonzero coprime positive integers $x,y$ and primes $p>3,~H_p(x,y)$ is not divisible by $p^2$}. Here we note that there is no condition on $x$ and $y$ except $(x,y)=1$.

\noindent
{\bf Proposition1:} {\it The first case of FLT is true for primes $p \equiv -1 (\bmod~6)$}. We will now show how Conjecture1 can be used to prove Proposition1. Assume the first case of FLT is false, then from Lemma2 $p^2 \mid G_p(x,y)$. But from Lemma3 we see that $p \nmid f_{x,y}$. So $p^2 \mid H_p(x,y)$ contradicting Conjecture1.

\noindent
{\bf Pollaczek's result [1]:} {\it If the non-zero integers $x,y,z$ satisfy equation $x^p+y^p=z^p,(x,y,z)=1 $, for prime $p>3$, then $x^2+xy+y^2$ is non-divisible by $p$}. 

\noindent
{\bf Proposition2:} {\it The first case of FLT is true for any prime $p>3$}. Again we will show how Conjecture1 can be used to prove Proposition2. As we have already seen if the first case of FLT is false then $p^2 \mid G_p(x,y)$. But from Pollaczek's result we see that $p \nmid f_{x,y}$. So $p^2 \mid H_p(x,y)$ contradicting Conjecture1.

\noindent
We propose a weaker {\it conjecture} which can be easily verified without altering the propositions.

\noindent
{\bf Conjecture2:} {\it For any prime $p>3$ assume that there exist pairwise relatively prime nonzero integers $x,y,z,x^p+y^p=z^p$. Then $H_p(x,y)$ is not divisible by $p^2$}.

\noindent
Since $H_5(x,y)=H_7(x,y)=1$, the conjectures are true for primes $5$ and $7$. For primes $p>7$ using a simple computer programme in Mathematica 4.0 the {\it conjectures} were being verified. Since $F_p(x,y)$ and $xy(x+y)f^\epsilon_{x,y}$ are symmetric polynomials, ${H_p(x,y)= \frac {F_p(x,y)}{pxy(x+y)f^\epsilon_{x,y}}}$ is also a symmetric polynomial. So $(y^{-1})^{p-2 \epsilon-3}~H_p(x,y) \equiv H_p(X,1) (\bmod~p~\hbox{or}~p^2)$, where $X \equiv y^{-1}x (\bmod~p~\hbox{or}~p^2)$ and $y^{-1}y \equiv 1 (\bmod~p~\hbox{or}~p^2)$. Due to the property $H_p(X,1)=H_p(-X-1,1)$ the number of residues, of $p$ or $p^2$, for the verification of the conjectures are reduced to half. For most primes $p$, we found, $p \nmid H_p(X,1)$ and the number of residues required for both conjectures is $(p-1)/2$. But when $p \mid H_p(X,1)$ verification of Conjecture1 takes longer than Conjecture2, since the number of residues of $p^2$ for Conjecture1 is $p(p-1)/2$ and for Conjecture2, because of condition from (12), $X^{p-1} \equiv 1 (\bmod~p^2)$, it is only $(p-1)/2$. With the computer Conjecture1 was verified to be true for primes $\leq 2437$ and Conjecture2 for primes $\leq 100003$. From Corollary1, the first case of FLT is true for prime $3$. So with Conjecture2 and Pollaczek's result, {\bf the first case of FLT is true for all primes $\leq 100003$}. We conclude this paper with a {\it corollary}.

\noindent
{\bf Corollary3:} {\it For any prime $p$, $2^{p-1}-1$ is not divisible by $p^3$.} This corollary can be easily checked for primes $2$ and $3$. For primes $p>3$, we have from Conjecture1 or Conjecture2 $p^2 \nmid H_p(x,y)$ even for $x=y=1$. But $H_p(x,y)= \frac {F_p(x,y)}{pxy(x+y)f^\epsilon_{x,y}}$ and substituting $x=y=1$ we get, $p^2 \nmid \frac{2^{p-1}-1}{3^\epsilon~p} \Rightarrow p^3 \nmid 2^{p-1}-1.~~~~~~~~~~~~~~~~~~~~~~~~~~~~~~~~~~~~~~~~~~~~~~~~~~~~~~~~~~~~~~~~~~~~~~~~~~~~~~~~~~~~~~~~\Box$

\noindent
{\bf Acknowledgement.} I am very thankful to M.A. Prasad for the discussions which helped in improvement of this paper. His suggestions in the computer programme reduced the run time considerbly. I am highly indebted to him for his valuable guidence and encouragment. I thank V. Nandagopal, School of Mathematics, Tata Institute of Fundamental Research, Mumbai, for helping in executing the programme.  

\vspace{0.3cm}
\noindent
{\bf REFERENCES}
\begin{enumerate}
\item Pollaczek, F., {\it $\ddot Uber~den~gro \beta en~Fermatschen~Satz$}, Sitzungsber. Akad. Wiss. Wein. Math.-Natur. Kl. IIa, {\bf 126} (1917), 45-59.
\end{enumerate}
\end{document}